# GENERALIZED SCORE TEST OF HOMOGENEITY FOR MIXED EFFECTS MODELS

By Hongtu Zhu and Heping Zhang[1]

*Columbia University and Yale University*

Many important problems in psychology and biomedical studies require testing for overdispersion, correlation and heterogeneity in mixed effects and latent variable models, and score tests are particularly useful for this purpose. But the existing testing procedures depend on restrictive assumptions. In this paper we propose a class of test statistics based on a general mixed effects model to test the homogeneity hypothesis that all of the variance components are zero. Under some mild conditions, not only do we derive asymptotic distributions of the test statistics, but also propose a resampling procedure for approximating their asymptotic distributions conditional on the observed data. To overcome the technical challenge, we establish an invariance principle for random quadratic forms indexed by a parameter. A simulation study is conducted to investigate the empirical performance of the test statistics. A real data set is analyzed to illustrate the application of our theoretical results.

**1. Introduction.** Mixed effects and latent variable models provide an attractive framework to accommodate correlated data. For example, structure equation models and generalized linear mixed models (GLMMs) are commonly used in behavioral, educational and social sciences (e.g., [2, 3]). A fundamental question in mixed effects or latent variable models is whether or not the inclusion of the random effects or latent variables is necessary. Many authors have examined this important issue using score test statistics in the framework of the GLMMs; see [8, 16, 21, 22, 32], for example. However, those authors did not fully exploit the general correlation structure of the random effects (or latent variables).

Received January 2004; revised June 2005.

[1]Supported in part by National Institute on Drug Abuse Grants DA12468, DA016750 and DA017713 and NSF Grant SES-05-50988.

Suppose that we observe data from $n$ units and within the $i$th unit we have $m_i$ measurements, $i = 1, \ldots, n$. This is a typical data structure in longitudinal and family studies that are popular in social and biomedical studies. In longitudinal studies, the unit is usually a person or an animal. In family studies, the unit is generally a family. In addition to the following two examples, other examples can be found in [38].

EXAMPLE 1 (Segregation analysis of ordinal traits). To study the genetic inheritance pattern of many health conditions such as cancer and psychiatric disorders, Zhang, Feng and Zhu [34] proposed a general framework for conducting complex segregation analysis of ordinal traits based on the latent variable model of Zhang and Merikangas [35]. Let $Y_i = (y_{i,1}, \ldots, y_{i,m_i})^T$ be a vector of traits and $X_i$ the covariates from the $i$th family, $i = 1, \ldots, n$. Without loss of generality, suppose that $y_{i,j}$ assumes ordinal values 0, 1 or 2. To model the potential familial correlation, they introduced a latent variable vector $\mathbf{v}_i$ for each family to represent common unmeasured environmental and genetic factors shared by family members. Conditional on the $\{\mathbf{v}_i\}$, the $y_{i,j}$'s are assumed to be independent and follow the proportional odds logistic model given by

$$\text{logit } P\{y_{i,j} = 0 | \mathbf{v}_i\} = \mathbf{x}_{i,j}^T \beta + \alpha_0 + b_{i,j},$$
$$(1)$$
$$\text{logit } P\{y_{i,j} \leq 1 | \mathbf{v}_i\} = \mathbf{x}_{i,j}^T \beta + \alpha_1 + b_{i,j},$$

where $\alpha_0 \leq \alpha_1$, $b_{i,j}$ depends on $\{\mathbf{v}_i\}$ and $\mathbf{x}_{i,j}$ is a covariate vector in the design matrix $X_i = (\mathbf{x}_{i,1}^T, \ldots, \mathbf{x}_{i,m_i}^T)$ $(m_i \times q_1)$ from the $j$th member in the $i$th family. An important objective in collecting family data is to test familial aggregation and inheritance, which can be achieved by testing $\text{var}[b_{i,j}] = 0$ for all $i$ and $j$.

EXAMPLE 2 (Generalized linear mixed effects model). Consider a data set that is composed of a response $y_{i,j}$, covariate vectors $\mathbf{x}_{i,j}$ and $\mathbf{z}_{i,j}$ for observations $j = 1, \ldots, m_i$ within clusters $i = 1, \ldots, n$. We define the generalized linear mixed effects models as

$$(2) \qquad p(y_{i,j} | \mathbf{b}_i) = \exp[\phi\{y_{i,j}\theta_{i,j} - a(\theta_{i,j})\} + c(y_{i,j}, \phi)]$$

and $\mu_{i,j} = E(y_{i,j} | \mathbf{b}_i) = g(\mathbf{x}_{i,j}^T \beta + \mathbf{z}_{i,j}^T \mathbf{b}_i)$, where $a(\cdot)$, $c(\cdot)$ and $g(\cdot)$ are known continuously differentiable functions. The random coefficients $\mathbf{b}_i$'s $(q \times 1)$ are normally distributed such that $E[\mathbf{b}_i] = 0$ and $E[\mathbf{b}_i \mathbf{b}_i^T] = \Sigma$. Moreover, for $i \neq i'$, $\mathbf{b}_i$ and $\mathbf{b}_{i'}$ are independent of each other. The so-called homogeneity test is to test whether $\Sigma = \mathbf{0}$.

To summarize the two examples presented above, we consider the following mixed effects model. We use $(y_{i,j}, \mathbf{x}_{i,j}, \mathbf{z}_{i,j})$ to denote the $j$th observation in the $i$th cluster. The total number of observations is $N = \sum_{i=1}^{n} m_i$. Furthermore, we assume that for each $Y_i$, there exists an unobserved $q \times 1$ latent variable (or random effect) vector $\mathbf{b}_i$. Given $\{\mathbf{b}_i; i = 1, \ldots, n\}$, the components of $\{Y_i; i = 1, \ldots, n\}$ are independent random variables and have the joint probability density function

$$(3) \qquad p(Y_i | \mathbf{b}_i) = \prod_{j=1}^{m_i} p(y_{i,j} | \psi_{i,j}(\mathbf{b}_i; \beta, \gamma_{(1)}), \Phi),$$

where $\psi_{i,j}(\mathbf{b}_i; \beta, \gamma_{(1)}) = g(\mathbf{x}_{i,j}^T \beta; f_{i,j}(\mathbf{z}_{i,j}, \gamma_{(1)})^T \mathbf{b}_i)$ and $\Phi$ is a dispersion parameter vector. In addition, $g(\cdot)$ is a known link function and $f_{i,j}(\cdot)$ is a $q \times 1$ vector function, and $\beta$ and $\gamma_{(1)}$ are, respectively, $q_1 \times 1$ and $q_3 \times 1$ vectors. The unobserved random variables, $\mathbf{b}_i$, satisfy $E[\mathbf{b}_i] = 0$ and $E[\mathbf{b}_i \mathbf{b}_{i'}^T] = \Sigma_{i,i'}(\gamma)$, where $\gamma$ is a $q_2 \times 1$ vector. Model (3) also includes the factor analysis model and the random coefficient model, in which $f_{i,j}(\cdot, \cdot)$ may depend on unknown parameters. Hereafter, we include $\gamma_{(1)}$ in $\gamma$ for notational simplicity.

We are interested in testing the homogeneity hypotheses

$$(4) \quad H_0 : \Sigma_{i,i'}(\gamma) = \mathbf{0} \qquad \text{for all } i, i' \quad \text{vs.} \quad H_1 : \Sigma_{i,i'}(\gamma) \neq \mathbf{0} \qquad \text{for some } i, i'.$$

We generally conduct the omnibus testing in (4), because it is easy to control its type I error. If the null hypothesis is rejected, it is interesting to find out which components are nonzero. While the details warrant a separate investigation, the results presented here will be useful for testing that some parameters in (4) equal zero.

To test the homogeneity hypotheses in (4), we need to address the following four issues: (a) a convenient parameterization for the homogeneity test; (b) the construction of a score test statistic; (c) the asymptotic distribution of the score test statistic under the null hypothesis; and (d) the computation of the $p$-value from the asymptotic distribution.

The solution to the first issue on the parameterization lays the foundation for resolving the subsequent issues. Let us examine a simple case of Example 2 with $q = 2$. We write the covariance matrix of $\mathbf{b}_i$, $\Sigma$ as

$$(5) \quad \begin{pmatrix} \sigma_1^2 & \rho \sigma_1 \sigma_2 \\ \rho \sigma_1 \sigma_2 & \sigma_2^2 \end{pmatrix} = \sigma_T \begin{pmatrix} \cos^2(\gamma_1) & \gamma_2 \sin(\gamma_1) \cos(\gamma_1) \\ \gamma_2 \sin(\gamma_1) \cos(\gamma_1) & \sin^2(\gamma_1) \end{pmatrix},$$

where $\sigma_T = \sigma_1^2 + \sigma_2^2$ and $(\sigma_1^2/\sigma_T, \sigma_2^2/\sigma_T) = (\cos^2(\gamma_1), \sin^2(\gamma_1))$. We see that the null hypothesis in (4) is equivalent to $\Sigma(\gamma) = \mathbf{0}$, that is, $\sigma_T = 0$. The first and second derivatives of the log-likelihood function with respect to all parameters $\sigma_1, \sigma_2$, and $\rho$ in $\Sigma(\gamma)$ are not continuous when $\Sigma(\gamma) = 0$; however, they are continuous in $\sigma_T$ at $\sigma_T = 0$ [1]. In this simple case, we simply test $\sigma_T = 0$ and treat the other parameters as nuisance parameters.

When $q \geq 2$, we consider a lower triangular Cholesky decomposition of $\Sigma$, denoted by $\mathbf{L} = (\ell_{i,j})$, which satisfies $\ell_{i,i} \geq 0$ for all $i = 1, \ldots, q$ and $\ell_{i,j} = 0$ for $i < j$. Furthermore, we define $\mathbf{L} = \Lambda\Gamma$, where $\Lambda = \operatorname{diag}(\ell_{1,1}, \ldots, \ell_{q,q})/\sqrt{\sum_{i=1}^{q} \ell_{i,i}^2}$ and $\Gamma = (\gamma_{i,j})$ is a $q \times q$ lower triangular matrix with $\gamma_{i,i} = 1$ for $i = 1, \ldots, q$. Let $\sigma_T = \sum_{i=1}^{q} \ell_{i,i}^2$. Then $\Sigma$ can be written as $\Sigma = \sigma_T \Lambda\Gamma\Gamma^T\Lambda$. Thus, the null hypothesis in (4) is equivalent to $\sigma_T = 0$.

To our knowledge, there are no satisfactory solutions to the remaining three issues. For example, Chen, Chen and Kalbfleisch [5, 6], Chen and Chen [4], Crainiceanu and Ruppert [10] and Zhu and Zhang [37] derived the asymptotic or small sample distributions of the likelihood ratio statistics for some specific mixed effects models under restrictive conditions. Others considered score test statistics. Liang [21] and Commenges and Jacqmin-Gadda [8] considered the case when the random effect, $\mathbf{b}_i$, is scalar. Although Lin [22] and Hall and Præstgaard [16] considered a multidimensional $\mathbf{b}_i$, their $f_{i,j}(\cdot, \cdot)$ does not contain $\gamma_{(1)}$. In other words, the existing results do not cover our examples and the general model (3). Thus, it is imperative for us to develop score test statistics and establish the asymptotic theory under a more general framework.

## 2. Score test statistics of homogeneity.

### 2.1. Score test statistics. From now on, we write

$$
(6) \qquad \Sigma_{i,i'}(\gamma) = \sigma_T W_{i,i'}(\gamma) \qquad \text{for all } i, i' = 1, \ldots, n,
$$

where $\sigma_T$ is introduced above. Under the parameterization (6), we formally state the homogeneity hypotheses as

$$
(7) \qquad H_0: \sigma_T = 0 \quad \text{vs.} \quad H_1: \sigma_T > 0.
$$

Letting $\mathbf{u}_i = \sigma_T^{-1/2} \mathbf{b}_i$, we see that $E[\mathbf{u}_i] = 0$ and $E[\mathbf{u}_i \mathbf{u}_{i'}^T] = W_{i,i'}(\gamma)$. Thus, the log-likelihood function $\mathcal{L}_n(\sigma_T | \beta, \gamma, \Phi)$ is given by

$$
\log \left\{ \int \prod_{i=1}^{n} \prod_{j=1}^{m_i} p(y_{i,j} | \psi_{i,j}(\mathbf{x}_{i,j}^T \beta; f_{i,j}(\mathbf{z}_{i,j}, \gamma_{(1)})^T \mathbf{u}_i \sigma_T^{1/2}), \Phi) \, dF(\mathbf{u}_1, \ldots, \mathbf{u}_n | \gamma) \right\},
$$

where $F(\mathbf{u}_1, \ldots, \mathbf{u}_n | \gamma)$ is the distribution function of $(\mathbf{u}_1, \ldots, \mathbf{u}_n)$. Let $t_{i,j} = \sigma_T^{1/2} \eta_{i,j}$, where $\eta_{i,j} = f_{i,j}(\mathbf{z}_{i,j}, \gamma_{(1)})^T \mathbf{u}_i$. Similarly to Liang [21], we can show that the first-order right derivative of $\mathcal{L}_n(\sigma_T | \beta, \gamma, \Phi)$ at $\sigma_T = 0$, denoted by $T_S(\gamma | \beta, \Phi)$, is given by

$$
0.5 \int \left[ \left\{ \sum_{i=1}^{n} \sum_{j=1}^{m_i} \frac{\partial \log p(y_{i,j} | \psi_{i,j}(\mathbf{x}_{i,j}^T \beta; t_{i,j}))}{\partial t_{i,j}}(0) \eta_{i,j} \right\}^2 \right.
$$
$$
\left. + \sum_{i=1}^{n} \sum_{j=1}^{m_i} [\eta_{i,j}]^2 \left\{ \frac{\partial^2 \log p(y_{i,j} | \psi_{i,j}(\mathbf{x}_{i,j}^T \beta; t_{i,j}))}{\partial t_{i,j}^2}(0) \right\} \right] dF(\mathbf{u}_1, \ldots, \mathbf{u}_n | \gamma);
$$

see [38] for a detailed derivation. We will describe later how to estimate $(\beta, \Phi)$, but for the time being, let us treat them as if they were known and not include them as parameters to simplify the notation. That is, let $T_S(\gamma) = T_S(\gamma | \beta, \Phi)$. If $\gamma$ is actually absent in all of the $W_{i,i'}(\gamma)$'s, $T_S(\gamma)$ is a score test statistic identical to that proposed by Liang [21] and Commenges and Jacqmin-Gadda [8]. In general, however, $T_S(\gamma)$ is not really a score statistic due to the presence of $\gamma$.

Let $b_{K,K'}(\gamma)$ be $f_{i,j}(\mathbf{z}_{i,j}, \gamma_{(1)})^T W_{i,i'}(\gamma) f_{i',j'}(\mathbf{z}_{i',j'}, \gamma_{(1)})$ and $B(\gamma) = (b_{K,K'}(\gamma))$ be an $N \times N$ matrix, where $K = (i, j)$ and $K' = (i', j')$. With this notation, $T_S(\gamma)$ can be decomposed into two terms,

$$(8) \qquad T_S(\gamma) = \mathbf{U}^T B(\gamma) \mathbf{U} - \mathrm{tr}[\mathbf{V} B(\gamma)],$$

where $\mathbf{U}$ and $\mathbf{V}$ are, respectively, an $N \times 1$ vector and an $N \times N$ matrix. Let $U_{i,j}$ and $V_{i,j}$ be the limits of $\partial \log p(y_{i,j} | \psi_{i,j}(\mathbf{x}_{i,j}^T \beta; t_{i,j})) / \partial t_{i,j}$ and $-\partial^2 \log p(y_{i,j} | \psi_{i,j}(\mathbf{x}_{i,j}^T \beta; t_{i,j})) / \partial t_{i,j}^2$, respectively, as $t_{i,j} \to 0$. The $K$th element of $\mathbf{U}$ is $U_{i,j}$, and $\mathbf{V}$ is a diagonal matrix with $K$th element $V_{i,j}$.

Following Commenges and Jacqmin-Gadda [8], we can decompose $T_S(\gamma)$ into two terms,

$$T_S(\gamma) = T_P(\gamma) + T_O(\gamma), \qquad T_O(\gamma) = \sum_K b_{K,K}(\gamma)(U_K^2 - V_K),$$

$$(9)$$

$$T_P(\gamma) = \mathbf{U}^T \{B(\gamma) - \mathrm{diag}[B(\gamma)]\} \mathbf{U} = \sum_{K \neq K'} b_{K,K'}(\gamma) U_K U_{K'},$$

where $\mathrm{diag}[B(\gamma)]$ is the $N \times N$ diagonal matrix of $B(\gamma)$. The first term $T_P(\gamma)$ is called a *pairwise correlation term* and the second term $T_O(\gamma)$ is an *overdispersion term*. Under the null hypothesis $H_0$, we have

$$E[T_P(\gamma)] = E[T_O(\gamma')] = E[T_P(\gamma) T_O(\gamma')] = 0,$$

$$E[T_S(\gamma) T_S(\gamma')] = E[T_P(\gamma) T_P(\gamma')] + E[T_O(\gamma) T_O(\gamma')],$$

$$E[T_P(\gamma) T_P(\gamma')] = 2 \sum_{K \neq K'} b_{K,K'}(\gamma') b_{K,K'}(\gamma) E U_K^2 E U_{K'}^2,$$

$$E[T_O(\gamma) T_O(\gamma')] = \sum_K b_{K,K}(\gamma) b_{K,K}(\gamma') [E U_K^4 + E V_K^2 - 2 E(U_K^2 V_K)].$$

We construct three score test statistics in the following. We first define

$$X_P(\gamma) = \frac{T_P(\gamma)}{\sqrt{I_{TP}(\gamma)}},$$

$$(10)$$

$$X_O(\gamma) = \frac{T_O(\gamma)}{\sqrt{I_{TO}(\gamma)}} \quad \text{and} \quad X_S(\gamma) = \frac{T_S(\gamma)}{\sqrt{I_{TS}(\gamma)}},$$

where $I_{TO}(\gamma)$, $I_{TS}(\gamma)$ and $I_{TP}(\gamma)$ are the variances of $T_O(\gamma)$, $T_S(\gamma)$ and $T_P(\gamma)$, respectively. However, we need to estimate $\xi = (\beta, \Phi)$ in $X_P(\gamma)$, $X_O(\gamma)$ and $X_S(\gamma)$ for testing and replace $\xi$ by its estimator $\hat{\xi}$. Let $\hat{U}_K$ and $\hat{V}_K$ denote the values of $U_K$ and $V_K$ evaluated at $\hat{\xi}$, respectively, which gives $\hat{T}_O(\gamma) = \sum_K b_{K,K}(\gamma)(\hat{U}_K^2 - \hat{V}_K)$, $\hat{T}_P(\gamma) = \sum_{K \neq K'} b_{K,K'}(\gamma)\hat{U}_K \hat{U}_{K'}$ and $\hat{T}_S(\gamma) = \hat{T}_P(\gamma) + \hat{T}_O(\gamma)$. We introduce

$$
\begin{aligned}
(11) \qquad & \hat{X}_P(\gamma) = \frac{\hat{T}_P(\gamma)}{\sqrt{I_{EP}(\gamma)}}, \\
& \hat{X}_O(\gamma) = \frac{\hat{T}_O(\gamma)}{\sqrt{I_{EO}(\gamma)}} \quad \text{and} \quad \hat{X}_S(\gamma) = \frac{\hat{T}_S(\gamma)}{\sqrt{I_{ES}(\gamma)}},
\end{aligned}
$$

where $I_{EO}(\gamma)$, $I_{ES}(\gamma)$ and $I_{EP}(\gamma)$ are the asymptotic variances of $\hat{T}_O(\gamma)$, $\hat{T}_S(\gamma)$ and $\hat{T}_P(\gamma)$, respectively, with $\xi$ evaluated at $\hat{\xi}$. Assume that

$$
(12) \qquad N^{1/2}(\hat{\xi} - \xi_*) = N^{-1/2} \sum_K F_K + o_p(1),
$$

where $\xi_*$ is the true value of $\xi$ and $F_K$ is a random function of $(y_K, \mathbf{x}_K, \mathbf{z}_K)$. In addition, $F_K$ and $F_{K'}$ are independent of each other for $K \neq K'$. With these preparations, we can show that $\hat{T}_O(\gamma) = \sum_K b_{K,K}(\gamma)(\hat{U}_K^2 - \hat{V}_K) \approx \sum_K \{b_{K,K}(\gamma)(U_K^2 - V_K) - J_N(\gamma)^T F_K\}$ under some mild conditions, where $J_N(\gamma) = E[-N^{-1} \partial_\xi T_O(\gamma)]$. Furthermore, $I_{EO}(\gamma)$ can be approximated by $I_{TO}(\gamma) - 2\sum_K b_{K,K}(\gamma)E[(U_K^2 - V_K)F_K^T]J_N(\gamma) + J_N(\gamma)^T \sum_K E(F_K F_K^T)J_N(\gamma)$.

Because of the one-sided constraint $\sigma_T \geq 0$, we consider $\hat{X}_P(\gamma)\mathbf{1}(\hat{X}_P(\gamma) \geq 0)$, $\hat{X}_O(\gamma)\mathbf{1}(\hat{X}_O(\gamma) \geq 0)$ and $\hat{X}_S(\gamma)\mathbf{1}(\hat{X}_S(\gamma) \geq 0)$, where $\mathbf{1}(A)$ is the indicator function of the event $A$. Furthermore, to remove the unknown $\gamma$, we introduce the maximum statistics defined by $S_O = \sup_\gamma \{\hat{X}_O(\gamma)^2 \mathbf{1}(\hat{X}_O(\gamma) \geq 0)\}$, $S_P = \sup_\gamma \{\hat{X}_P(\gamma)^2 \mathbf{1}(\hat{X}_P(\gamma) \geq 0)\}$ and $S_S = \sup_\gamma \{\hat{X}_S(\gamma)^2 \mathbf{1}(\hat{X}_S(\gamma) \geq 0)\}$. In practice, the null hypothesis is rejected if any of these three statistics $\{S_O, S_P, S_S\}$ has a large absolute value.

As a common practice, the foregoing use of the maximum of the score test statistics is based on power considerations (see, e.g., [14]). Because $\gamma$ is identifiable under the alternative hypothesis only, the maximization over $\gamma$ takes effect under the alternative hypothesis, as for the likelihood ratio test (LRT). We show in [38] that $S_S$ yields an efficient test statistic because it recovers information from the likelihood under the alternative hypothesis. Furthermore, we show that the score test statistic proposed here is asymptotically equivalent to the LRT for testing the homogeneity of random effects; see Theorems S.1 and S.2 in [38].

By now we have defined three score statistics for testing homogeneity under mixed effects models, but we will discuss their asymptotic null distributions in Section 3. Similar to Lin's [22] method, an important feature of our score statistics is that we only need to specify the first and second moments of the latent variables in (8) for the distribution function $F(\mathbf{b}_1, \ldots, \mathbf{b}_n; \gamma)$. Thus, the test statistics are expected to be robust with respect to the distribution of the random effects. In addition, our test statistics allow a general covariance structure of the latent variables, and $f_{i,j}(\cdot, \cdot)$ may depend on unknown parameters.

As we know, the optimality of a test depends on its power. To compare the power of $S_Q, S_P$ and $S_S$ with that of Lin [22], we consider sequences of local alternatives to $\sigma_T = 0$. The asymptotic local power for $S_S$ follows from Theorem 2 in [14]; see Theorems S.1–S.4 in [38] for details. Empirically, simulations in Section 4 will demonstrate that the score statistic $S_S$ proposed here is more powerful than the score statistic proposed by Lin [22] (see Tables 1 and 2).

2.2. *A resampling procedure.* To assess the power of the three test statistics $\{S_S, S_P, S_O\}$, we need to obtain empirical distributions for the score statistics in lieu of their theoretical distributions. What follow are the four key steps in generating the stochastic processes that have the same asymptotic distributions as the test statistics.

*Step* 1. We generate i.i.d. random samples, $\{v_{K,K'}^{(r)} : K, K' = (i,j), j = 1, \ldots, m_i, i = 1, \ldots, n\}$, from $N(0,1)$. Here, the superscript $(r)$ represents a replication number.

*Step* 2. We calculate $\hat{T}_S^{(r)}(\gamma) = \hat{T}_P^{(r)}(\gamma) + \hat{T}_O^{(r)}(\gamma)$ and

$$
\begin{aligned}
\hat{T}_P^{(r)}(\gamma) &= \sqrt{2} \sum_{K \neq K'} b_{K,K'}(\gamma) \hat{U}_K \hat{U}_{K'} v_{K,K'}^{(r)}, \\
\hat{T}_O^{(r)}(\gamma) &= \sum_K v_{K,K}^{(r)} \{ b_{K,K}(\gamma)(\hat{U}_K^2 - \hat{V}_K) - J_N(\gamma)^T \hat{F}_K \},
\end{aligned}
\tag{13}
$$

where $\hat{F}_K$ is an estimator of $F_K$ evaluated at $\hat{\xi}$. Then, we can calculate

$$
\hat{X}_S^{(r)}(\gamma) = \frac{\hat{T}_S^{(r)}(\gamma)}{\sqrt{I_{ES}(\gamma)}},
$$

$$
\hat{X}_P^{(r)}(\gamma) = \frac{\hat{T}_P^{(r)}(\gamma)}{\sqrt{I_{EP}(\gamma)}} \quad \text{and} \quad \hat{X}_O^{(r)}(\gamma) = \frac{\hat{T}_O^{(r)}(\gamma)}{\sqrt{I_{EO}(\gamma)}}.
$$

It is important to note that conditional on the observed data, $\hat{X}_S^{(r)}(\gamma)$, $\hat{X}_P^{(r)}(\gamma)$ and $\hat{X}_O^{(r)}(\gamma)$ converge weakly to the three Gaussian processes described in Theorem 2 as $N \to \infty$ (see Section 3). This can be shown using the conditional functional central limit theorem; see Section 3 for details.

*Step* 3. We calculate the three test statistics

$$
S_S^{(r)} = \sup_{\gamma \in \Gamma} \{ \hat{X}_S^{(r)}(\gamma)^2 \mathbf{1}(\hat{X}_S^{(r)}(\gamma) \geq 0) \}, \qquad S_P^{(r)} = \sup_{\gamma \in \Gamma} \{ \hat{X}_P^{(r)}(\gamma)^2 \mathbf{1}(\hat{X}_P^{(r)}(\gamma) \geq 0) \}
$$

and

$$S_O^{(r)} = \sup_{\gamma \in \Gamma} \{\hat{X}_O^{(r)}(\gamma)^2 \mathbf{1}(\hat{X}_O^{(r)}(\gamma) \geq 0)\}.$$

*Step* 4. We repeat the above three steps $r_0$ times and obtain three realizations: $\{S_S^{(r)} : r = 1, \ldots, r_0\}$, $\{S_P^{(r)} : g = 1, \ldots, r_0\}$ and $\{S_O^{(r)} : r = 1, \ldots, r_0\}$. It can be shown that the empirical distribution of $S_S^{(r)}$ converges to the asymptotic distribution of $S_S$. Similarly, $S_O^{(r)}$ and $S_P^{(r)}$ converge to the asymptotic distributions of $S_O$ and $S_P$, respectively. Therefore, the empirical distributions of these three realizations form the basis for calculating the significance level and power of the tests.

2.3. *Example* 2 (Continued). Let us revisit Example 2 to illustrate how our test statistics can be applied. Using the parameterization in Section 1, we see that $\Sigma(\gamma) = \sigma_T W(\gamma)$, and $W_{i,i'}(\gamma)$ equals $W(\gamma)$ for $i = i'$ and is zero otherwise. In this case, $\mathbf{u}_i = \sigma_T^{-1/2}\mathbf{b}_i$, $\eta_{i,j} = \mathbf{z}_{i,j}^T\mathbf{u}_i$, $f_{i,j}(\mathbf{z}_{i,j}, \gamma_{(1)}) = \mathbf{z}_{i,j}$ and $t_{i,j} = \sigma_T^{1/2}\eta_{i,j} = \mathbf{z}_{i,j}^T\mathbf{b}_i$. Moreover, we define $\theta_{i,j}(t_{i,j}) = k(\mathbf{x}_{i,j}^T\beta + t_{i,j})$ and $\mu_{i,j}(t_{i,j}) = g(\mathbf{x}_{i,j}^T\beta + t_{i,j})$ to emphasize the fact that they depend on $t_{i,j}$ explicitly. After some calculations, for model (2) we have $U_{i,j} = \phi e_{i,j}\dot{k}(\mathbf{x}_{i,j}^T\beta + t_{i,j})|_{t_{i,j}=0}$ and

$$V_{i,j} = \{\phi\ddot{a}(\theta_{i,j}(0))\dot{k}(\mathbf{x}_{i,j}^T\beta + t_{i,j})^2 - \phi e_{i,j}\ddot{k}(\mathbf{x}_{i,j}^T\beta + t_{i,j})\}|_{t_{i,j}=0},$$

where the dots denote differentiation, $\theta_{i,j}(0) = k(\mathbf{x}_{i,j}^T\beta + t_{i,j})|_{t_{i,j}=0}$ and $e_{i,j} = y_{i,j} - \mu_{i,j}(0) = y_{i,j} - \mu_{i,j}(t_{i,j})|_{t_{i,j}=0}$. In addition, $b_{K,K'}(\gamma) = \mathbf{z}_{i,j}^T W(\gamma)\mathbf{z}_{i,j'}$ for $i = i'$ and is zero otherwise. From now on, we use $\dot{k}_{i,j}$ to denote $\dot{k}(\mathbf{x}_{i,j}^T\beta + t_{i,j})$ evaluated at $t_{i,j} = 0$ and likewise for $\ddot{k}_{i,j}$. Thus, we have $T_S(\gamma) = T_O(\gamma) + T_P(\gamma)$, $T_P(\gamma) = \sum_{i=1}^n \sum_{j \neq j'} \mathbf{z}_{i,j}^T W(\gamma)\mathbf{z}_{i,j'}\phi^2 e_{i,j}e_{i,j'}\dot{k}_{i,j}\dot{k}_{i,j'}$ and

$$T_O(\gamma) = \sum_{i=1}^n \sum_{j=1}^{m_i} \mathbf{z}_{i,j}^T W(\gamma)\mathbf{z}_{i,j}[\phi^2 e_{i,j}^2 \dot{k}_{i,j}^2 + \phi e_{i,j}\ddot{k}_{i,j} - \phi\ddot{a}(\theta_{i,j}(0))\dot{k}_{i,j}^2].$$

Under the null hypothesis $H_0$, we use the first four central moments of $y_{i,j}$ of the exponential-family distributions [19, 33] to get

$$EU_{i,j}^2 = \phi\ddot{a}(\theta_{i,j}(0))\dot{k}_{i,j}^2,$$
$$EU_{i,j}^4 = \dot{k}_{i,j}^4\{3\phi^2\ddot{a}(\theta_{i,j}(0))^2 + \phi a^{(4)}(\theta_{i,j}(0))\},$$
$$EV_{i,j}^2 = \phi\ddot{a}(\theta_{i,j}(0))\ddot{k}_{i,j}^2 + \phi^2\ddot{a}(\theta_{i,j}(0))^2\dot{k}_{i,j}^4 \quad \text{and}$$
$$E(U_{i,j}^2 V_{i,j}) = \phi^2 a^{(3)}(\theta_{i,j}(0))\dot{k}_{i,j}^2\ddot{k}_{i,j} + \phi^2\ddot{a}(\theta_{i,j}(0))^2\dot{k}_{i,j}^4.$$

Thus, we can have $I_{TS}(\gamma) = I_{TO}(\gamma) + I_{TP}(\gamma)$ and

$$I_{TP}(\gamma) = 2\sum_{i=1}^{n}\sum_{j,j'=1, j\neq j'}^{m_i}\{\mathbf{z}_{i,j}^T W(\gamma)\mathbf{z}_{i,j'}\}^2\phi^2\ddot{a}(\theta_{i,j}(0))\ddot{a}(\theta_{i,j'}(0))\dot{k}_{i,j}^2\dot{k}_{i,j'}^2,$$

$$I_{TO}(\gamma) = \sum_{i=1}^{n}\sum_{j=1}^{m_i}\{\mathbf{z}_{i,j}^T W(\gamma)\mathbf{z}_{i,j}\}^2[EU_{i,j}^4 + EV_{i,j}^2 - 2E(U_{i,j}^2 V_{i,j})].$$

As discussed above, we need to replace $\beta$ and $\phi$ by their estimates under $H_0$. The maximum likelihood estimate, $\hat{\beta}$, of $\beta$ satisfies

$$\hat{\beta} - \beta_* = \left\{\sum_{i=1}^{n}\sum_{j=1}^{m_i}\ddot{a}(\theta_{i,j}(0))\dot{k}_{i,j}^2\mathbf{x}_{i,j}^T\mathbf{x}_{i,j}\right\}^{-1}\sum_{i=1}^{n}\sum_{j=1}^{m_i}\dot{k}_{i,j}e_{i,j}\mathbf{x}_{i,j} + o_p(N^{-1/2}),$$

which gives $F_K$ for each $K = (i, j)$. Moreover, we can calculate that

$$J_N(\gamma) = N^{-1}\sum_{i=1}^{n}\sum_{j=1}^{m_i}\mathbf{z}_{i,j}^T W(\gamma)\mathbf{z}_{i,j}\phi\{\ddot{a}(\theta_{i,j}(0))\dot{k}_{i,j}\ddot{k}_{i,j} + a^{(3)}(\theta_{i,j}(0))\dot{k}_{i,j}^3\}\mathbf{x}_{i,j},$$

and $I_{EO}(\gamma) = I_{TO}(\gamma) - J_N(\gamma)^T\{\sum_{i=1}^{n}\sum_{j=1}^{m_i}\ddot{a}(\theta_{i,j}(0))\dot{k}_{i,j}^2\mathbf{x}_{i,j}^T\mathbf{x}_{i,j}\}^{-1}J_N(\gamma)$, in which $\beta$ is replaced by $\hat{\beta}$. The strategy to deal with the unknown $\phi$ is similar.

**3. Asymptotic null distribution of score test statistics.** In this section, we study the asymptotic properties of $\{X_P(\gamma), \hat{X}_P(\gamma), \hat{X}_P^{(r)}(\gamma)\}$ under the null hypothesis $H_0$. Note that the asymptotic distributions of $X_O(\gamma)$, $\hat{X}_O(\gamma)$ and $\hat{X}_O^{(r)}(\gamma)$ have been widely discussed in the literature [25]. The asymptotic distribution of $\hat{X}_S^{(r)}(\gamma)$ follows from that of $\hat{X}_P^{(r)}(\gamma)$ and $\hat{X}_O^{(r)}(\gamma)$. We refer to [38] for details on how to apply our asymptotic results in some specific examples.

3.1. *Asymptotic null distribution.* We denote $\Rightarrow$ for weak convergence of a sequence of stochastic processes indexed by $\gamma \in \Gamma$, where the parametric space $\Gamma$ is a uniformly bounded convex compact subset of $R^{q_2}$. In addition, the uniform metric is used to define the weak convergence. Moreover, for a metric space $\{\mathcal{D}, d\}$, we consider $BL_1(\mathcal{D})$ to be the space of real-valued functions on $\mathcal{D}$ with Lipschitz norm bounded by 1, that is, for any $h \in BL_1(\mathcal{D})$, $\sup_{x\in D}|h(x)| \leq 1$ and $|h(x) - h(y)| \leq d(x, y)$. As discussed in [29], as $N \to \infty$, a stochastic process, $X_P^{(r)}$, weakly converges to $G_P$ on $\mathcal{D}$ if and only if $\sup_{h\in BL_1(D)}|Eh(X_P^{(r)}) - Eh(G_P)| \to 0$.

We have the following theorems, but we defer the proofs of all theorems as well as the assumptions to the Appendix.

Because $X_O(\gamma)$ can be regarded as the sum of independent but not identically distributed random variables indexed by $\gamma$, the asymptotic distribution of $X_O(\gamma)$ is a Gaussian process under some mild conditions by directly applying the functional central limit theorem (FCLT) [25, 29]. Furthermore, after examining the expressions for $X_P(\gamma)$ and $X_S(\gamma)$, we find that both of them are random quadratic forms indexed by $\gamma$. Although the asymptotic properties of random quadratic forms have been extensively studied in the literature (e.g., [12, 24]), those results are not applicable to $X_S(\gamma)$ and $X_P(\gamma)$, because these stochastic processes are indexed by $\gamma$. Thus, an invariance principle for the quadratic form process indexed by $\gamma$ needs to be developed; see detailed discussion in Section 3.2.

THEOREM 1. *Under conditions* (A1)–(A5) *in the Appendix and the null hypothesis $H_0$, $X_P(\cdot)$, $X_O(\cdot)$ and $X_S(\cdot)$ converge weakly to centered Gaussian processes as $N \to \infty$.*

Theorem 1 characterizes the asymptotic null distributions of the stochastic processes of interest and forms the foundation for constructing test statistics.

Let us understand the asymptotic properties of $\hat{X}_O(\gamma)$, $\hat{X}_P(\gamma)$ and $\hat{X}_S(\gamma)$. Because $\hat{X}_O(\gamma)$ is also asymptotically equivalent to the sum of independent random variables, it converges to a Gaussian process under some mild conditions (see Theorem 10.6 of [25]). For $\hat{X}_P(\gamma)$, under suitable conditions we can show that $\hat{X}_P(\gamma) = \sum_{K \neq K'} b_{K,K'}(\gamma) \hat{U}_K \hat{U}_{K'} / \sqrt{I_{EP}(\gamma)} = X_P(\gamma) + o_p(1)$ and $I_{EP}(\gamma) = I_{TP}(\gamma) + o_p(1)$. Thus, $\hat{X}_P(\gamma)$ and $X_P(\gamma)$ are asymptotically equivalent. The asymptotic distribution of $\hat{X}_S(\gamma)$ can be established by noting that it is a weighted sum of $\hat{X}_O(\gamma)$ and $\hat{X}_P(\gamma)$. To summarize our discussions, we have the following theorem.

THEOREM 2. *Under conditions* (B1)–(B8) *in the Appendix and the null hypothesis $H_0$, as $N \to \infty$, $\hat{X}_P(\cdot) \Rightarrow G_P(\cdot)$, $\hat{X}_O(\cdot) \Rightarrow G_O(\cdot)$ and $\hat{X}_S(\cdot) \Rightarrow G_S(\cdot)$, where $G_P$, $G_O$ and $G_S$ are three centered Gaussian processes.*

Theorem 2 delineates the asymptotic distributions of $\hat{X}_P(\gamma)$, $\hat{X}_O(\gamma)$ and $\hat{X}_S(\gamma)$. In the generalized linear models, $\hat{X}_O(\gamma)$ is the same as several tests for overdispersion [9]. In an example of a Bernoulli response variable, Jacqmin-Gadda and Commenges [17] show that $\hat{T}_S(\gamma)$ is identical to the pairwise correlation term.

To derive asymptotic null distributions of $S_O$, $S_P$ and $S_S$, we apply the continuous mapping theorem and have the following corollary.

COROLLARY. *Under the assumptions of Theorem 2, $S_O \xrightarrow{d} \sup_\gamma \{G_O(\gamma)^2 \times \mathbf{1}(G_O(\gamma) \geq 0)\}$, $S_P \xrightarrow{d} \sup_\gamma \{G_P(\gamma)^2 \mathbf{1}(G_P(\gamma) \geq 0)\}$ and $S_S \xrightarrow{d} \sup_\gamma \{G_S(\gamma)^2 \times \mathbf{1}(G_S(\gamma) \geq 0)\}$, where $\xrightarrow{d}$ represents convergence in distribution as $N \to \infty$.*

3.2. *Asymptotic distribution of a random quadratic form.* As noted above, to understand the asymptotic null distribution of $\{X_P(\gamma), \hat{X}_P(\gamma), \tilde{X}_P^{(r)}(\gamma)\}$, we need to investigate the asymptotic properties of the random quadratic forms indexed by $\gamma \in \Gamma$. For convenience, we will also use $K$ to index the integers from 1 to $N$ as well as the pairs $(i, j)$, because there is a one-to-one correspondence between $\{K = (i, j) : j = 1, \ldots, m_i, i = 1, \ldots, n\}$ and $\{K = 1, \ldots, N\}$. Consider the quadratic form without diagonal terms

$$(14) \qquad X_P(\gamma) = \sum_{K \neq K'} c_{K,K'}(\gamma) x_K x_{K'},$$

where $x_1, \ldots, x_N$ are a sequence of independent random variables such that $Ex_K = 0$ and $Ex_K^2 = 1$ for all $K = 1, \ldots, N$. Note that the $c_{K,K'}(\gamma)$'s may depend on $N$. We establish the asymptotic distribution of $X_P(\gamma)$ as follows.

THEOREM 3. *Under assumptions* (C1)–(C4) *in the* [Appendix]*,* $X_P \Rightarrow G_P$*, where* $G_P$ *is a centered Gaussian process with covariance matrix* $\rho_{(1)}(\gamma, \gamma')$*, as* $N \to \infty$*.*

Theorem [3] establishes an invariance principle for a random quadratic form indexed by $\gamma$; however, generalizing this result to a random quadratic form indexed by an arbitrary index set warrants further investigation. To simulate the asymptotic distribution of $X_P(\gamma)$, we consider the quadratic form

$$(15) \qquad X_P^{(r)}(\gamma) = \sqrt{2} \sum_{K \neq K'} c_{K,K'}(\gamma) x_K x_{K'} v_{K,K'}^{(r)},$$

where $\{v_{K,K'}^{(r)} : K, K' = 1, \ldots, N\}$ is a sequence of random variables defined in Step 1 above. Let $E_V$ denote the expectation taken with respect to all $v_{K,K'}^{(r)}$ conditional on the data.

THEOREM 4. *Assume that* (C2)–(C4) *in the* [Appendix] *are true and* (C1) *holds for* $p \geq 4$*. Then* $X_P^{(r)}(\cdot)$ *converges weakly to the same Gaussian process* $G_P(\cdot)$ *as* $N \to \infty$*; that is,* $X_P^{(r)}$ *is asymptotically measurable. In particular, as* $N \to \infty$*,* $\sup_{h \in BL_1(\ell^\infty(\Gamma))} |E_V h(X_P^{(r)}) - Eh(G_P)| \to 0$*, in probability.*

One important feature of Theorem [4] is that we can use $X_P^{(r)}$ to approximate the Gaussian process $G_P$. This theorem generalizes the resampling technique from the independent but nonidentically distributed framework [20] to the more general random quadratic setting. In particular, we propose a practical resampling technique to simulate the asymptotic distribution of $X_P(\gamma)$.

To consider the process $\hat{X}_P(\gamma)$, we introduce a sequence of independent random functions $\{U_1(s_1, \xi), \ldots, U_N(s_N, \xi)\}$ such that $U_K(s_K, \xi_*) = x_K$. Furthermore, let $X_P(\gamma, \xi - \xi_*) = \sum_{K \neq K'} c_{K,K'}(\gamma)U_K(s_K, \xi)U_{K'}(s_{K'}, \xi)$. We can easily see that $X_P(\gamma, \mathbf{0}) = X_P(\gamma)$ and $\hat{X}_P(\gamma) = X_P(\gamma, \hat{\xi} - \xi_*)$. In the following, we will prove that if $\hat{\xi} = \xi_* + O_p(N^{-1/2})$, then as $N \to \infty$,

$$(16) \qquad \hat{X}_P(\gamma) = X_P(\gamma, \hat{\xi} - \xi_*) = X_P(\gamma) + o_p(1).$$

Thus, the asymptotic distribution of $\hat{X}_P(\gamma)$ is the same as that of $X_P(\gamma)$ as described in Theorem 3. The asymptotic distribution of $\hat{X}_P(\gamma)$ in Theorem 2 follows directly from (16). A sufficient condition for (16) is that

$$(17) \qquad \sup_{\gamma \in \Gamma, \|\mathbf{h}\|_2 \leq M} |X_P(\gamma, \mathbf{h}N^{-1/2}) - X_P(\gamma, \mathbf{0})| = o_p(1)$$

holds for any given $M > 0$, where $\|\cdot\|_2$ is the Euclidean norm. The following theorem validates this sufficiency condition.

THEOREM 5. *Under assumptions* (C1)–(C8) *in the Appendix and* $\hat{\xi} = \xi_* + O_p(N^{-1/2})$, $\hat{X}_P(\gamma)$ *is asymptotically equivalent to* $X_P(\gamma)$ *as* $N \to \infty$. *In particular,* (17) *is true.*

Theorem 5 first gives the exact conditions to guarantee the asymptotic equivalence between $\hat{X}_P(\gamma)$ and $X_P(\gamma)$. Similarly, we can use $X_P^{(r)}(\gamma, \xi - \xi_*) = \sqrt{2} \sum_{K \neq K'} c_{K,K'}(\gamma) v_{K,K'}^{(r)} U_K(s_K, \xi) U_{K'}(s_{K'}, \xi)$ to approximate the asymptotic distribution of $\hat{X}_P(\gamma)$. In particular, $X_P^{(r)}(\gamma, \hat{\xi} - \xi_*)$ has the same form as $\hat{X}_P^{(r)}(\gamma)$ in Step 3 of Section 2.3. By using similar arguments to those in Theorem 5, we can prove that $\hat{X}_P^{(r)}(\gamma) = X_P^{(r)}(\gamma) + o_p(1)$. As shown in Theorem 4, $X_P^{(r)}$ converges to the process $G_P$ in distribution conditional on the data. Combining Theorems 4 and 5, we can conclude that $\hat{X}_P^{(r)}(\gamma)$ has desired properties, which leads to the following corollary.

COROLLARY. *Under assumptions* (C1)–(C8), $\hat{X}_P^{(r)}(\gamma)$ *is asymptotically equivalent to* $X_P^{(r)}(\gamma)$ *and* $\hat{X}_P^{(r)} \Rightarrow G_P$ *conditional on the data.*

**4. Simulation study and a real example.** There are two computational issues related to our test procedures in Section 2. First, we need to replace $\xi$ by $\hat{\xi}$. In the following, we choose $\hat{\xi}$ to be the maximum likelihood estimate obtained from the Newton–Raphson algorithm under the null hypothesis. Second, the computation for generating the three realizations, as required in Step 4 of the resampling process, is intensive. For instance, to generate $\{S_S^{(r)} : r = 1, \ldots, r_0\}$, each $S_S^{(r)}$ entails a maximization process because $S_S^{(r)} = \sup_{\gamma \in \Gamma} \{\hat{X}_S^{(r)}(\gamma)^2 \mathbf{1}(\hat{X}_S^{(r)}(\gamma) \geq 0)\}$. To ease the computational burden, we approximate $\Gamma$ by a grid $\Gamma_A$.

4.1. *A simulation study.* In this section, we use simulations to compare the performance of $S_P$, $S_O$ and $S_S$, and the test of Lin [22], denoted as $LS$.

The simulated data sets were drawn from two generalized linear mixed models: the logistic mixed model and the linear mixed model. We assume that the logistic mixed model has the form

$$(18) \quad \text{logit } P(y_{i,j} = 1 | \mathbf{b}_i) = 1.0 + 0.8x_{i,j1} + 0.5x_{i,j2} + (b_{i,1} + z_{i,j1}b_{i,2}).$$

The linear mixed model takes the form

$$(19) \quad y_{i,j} = 1.0 + x_{i,j1} + x_{i,j2} + (b_{i,1} + z_{i,j1}b_{i,2}) + \varepsilon_{i,j},$$

where $\varepsilon_{i,j}$ and the random effects $\mathbf{b}_i = (b_{i,1}, b_{i,2})^T$ are independent, and $\varepsilon_{i,j}$ follows a normal distribution with mean zero and variance $\phi$. The $x_{i,j1}$, $x_{i,j2}$ and $z_{i,j1}$ were simulated from a standard normal generator. The random effects $\mathbf{b}_i$ were simulated from a bivariate normal distribution with mean $(0,0)$ and the $2 \times 2$ covariance matrix $\sigma_1^2(1, \rho_1 \| \rho_1, \rho_2)$. Using the parameterization (5), $\Gamma$ is given by $(0, \pi] \times [-\delta_0, \delta_0]$, where $\delta_0$ is any scalar in $[0, 1)$. We used the grid $\Gamma_A = \{(i\pi/20, j/16) : i = 1, \ldots, 20; j = -15, \ldots, 15\}$. The size of the grid was based on computational feasibility and our empirical observation that it appears large enough to perform well. In the resampling procedure, $r_0$ was set to be 1,000.

For all score test statistics, we first compare the type I error under the null hypothesis and the power under the alternative hypotheses. For the logistic mixed model (18), we generated observations from the Bernoulli distribution $B(1, P(y_{i,j} = 1 | \mathbf{b}_i))$. We considered $n = 30$ and 50. Every unit contains 5 subjects ($m_i = 5$). We used correlations $\rho_1 = 0.5$ and $\rho_2 = 1.0$, and several different values of $\sigma_1$. For the linear mixed model (19), we demonstrate the gain of power by considering the correlation structure among the random effects. The simulated data set contains 40 ($n$) 4-subject units. To generate the random effects, we chose seven different values of $\sigma_1^2$ and two sets of $(\rho_1, \rho_2)$ given by $(0.5, 1)$ and $(-0.3, 0.2)$.

The results based on 10,000 replications are reported in Table 1. As expected, a larger sample size improves the power of detecting heterogeneity. The rejection rate under the null hypothesis is close to the nominal level of 0.05 for the score test statistics. Under the alternative hypothesis, $S_S$ is slightly and consistently more powerful than $S_P$. This is because $S_S$ accounts for the overdispersion due to the latent variables, which is tested by $S_O$. Table 1 also suggests that power is improved for all test statistics when a general correlation structure is assumed. It is quite remarkable that, even with relatively small sample sizes (30–50), the increment of power is still evident. Not surprisingly, the loss of power is more severe by ignoring the correlation when it is actually high. We should clarify that we used a general correlation structure to obtain the results in Table 1, instead of the underlying correlation. Also important, in all cases $S_S$ is more powerful than $LS$, and the difference is sometimes substantial. This observation is consistent with the fact that the likelihood ratio statistic under the constrained alternative is uniformly more powerful than that for the unconstrained case and that there is an asymptotic equivalence between the likelihood ratio statistic and the score statistic under both the constrained and unconstrained alternatives [27, 28, 38].

TABLE 1

*Comparison of type* I *error and the power for the score test statistics under models* (18) *and* (19) *at significance level* 0.05. *"Considered" and "ignored" represent including or excluding the correlation in the score test statistics*

| $\sigma_1^2$ | $LS$ | $\hat{S}_O$ | $\hat{S}_P$ | $\hat{S}_S$ | $LS$ | $\hat{S}_O$ | $\hat{S}_P$ | $\hat{S}_S$ |
|---|---|---|---|---|---|---|---|---|
| | | | logistic mixed model (18) and Bernoulli distribution | | | | | |
| | | $n = 30$ $(m_i = 5)$ | | | | $n = 50$ $(m_i = 5)$ | | |
| 0.0 | 0.049 | 0.048 | 0.054 | 0.060 | 0.052 | 0.056 | 0.048 | 0.054 |
| 0.3 | 0.216 | 0.060 | 0.282 | 0.292 | 0.312 | 0.076 | 0.426 | 0.446 |
| 0.6 | 0.484 | 0.064 | 0.600 | 0.608 | 0.681 | 0.061 | 0.760 | 0.780 |
| 0.8 | 0.626 | 0.060 | 0.734 | 0.734 | 0.840 | 0.070 | 0.908 | 0.906 |
| 1.2 | 0.828 | 0.050 | 0.884 | 0.886 | 0.969 | 0.092 | 0.972 | 0.972 |
| | | | linear mixed model (19) | | | | | |
| | | | $\Sigma = \sigma_1^2(1, 0.5\|0.5, 1)$ | | | | | |
| | | | $n = 40$ $(m_i = 4)$ | | | | | |
| | | ignored | | | | considered | | |
| 0.00 | 0.047 | 0.022 | 0.054 | 0.053 | 0.058 | 0.030 | 0.051 | 0.055 |
| 0.05 | 0.186 | 0.063 | 0.203 | 0.233 | 0.232 | 0.054 | 0.213 | 0.265 |
| 0.10 | 0.433 | 0.116 | 0.389 | 0.477 | 0.483 | 0.106 | 0.414 | 0.538 |
| 0.15 | 0.639 | 0.172 | 0.609 | 0.700 | 0.689 | 0.158 | 0.632 | 0.719 |
| 0.20 | 0.788 | 0.225 | 0.757 | 0.819 | 0.812 | 0.233 | 0.771 | 0.845 |
| 0.25 | 0.873 | 0.282 | 0.845 | 0.894 | 0.898 | 0.284 | 0.847 | 0.915 |
| 0.30 | 0.923 | 0.332 | 0.902 | 0.940 | 0.944 | 0.341 | 0.905 | 0.952 |
| | | | $\Sigma = \sigma_1^2(1, -0.3\|-0.3, 0.2)$ | | | | | |
| | | | $n = 40$ $(m_i = 4)$ | | | | | |
| | | ignored | | | | considered | | |
| 0.05 | 0.113 | 0.026 | 0.140 | 0.153 | 0.121 | 0.040 | 0.146 | 0.165 |
| 0.10 | 0.217 | 0.031 | 0.289 | 0.309 | 0.267 | 0.039 | 0.313 | 0.332 |
| 0.15 | 0.397 | 0.038 | 0.475 | 0.496 | 0.452 | 0.058 | 0.479 | 0.518 |
| 0.20 | 0.586 | 0.043 | 0.638 | 0.646 | 0.606 | 0.081 | 0.639 | 0.675 |
| 0.25 | 0.605 | 0.051 | 0.751 | 0.770 | 0.736 | 0.094 | 0.762 | 0.795 |
| 0.30 | 0.797 | 0.053 | 0.848 | 0.849 | 0.828 | 0.114 | 0.853 | 0.880 |

We now examine the type I error of all score tests in the case of nuisance overdispersion. For model (18), the binomial distribution, $B(5, P(y_{i,j} = 1|\mathbf{b}_i))$, was included in this study to introduce large overdispersion. We perturbed model (18) with random intercepts, which are independent of each other and subject specific. Specifically, we added $\sigma_2^2 v_{i,j}$ to the constant intercept 1 in model (18), where $v_{i,j}$ was generated from $N(0,1)$. For model (18), we simulated random errors $\varepsilon_{i,j}$ from $N(0, \phi \exp(v_{i,j}\sigma_2^2))$. We set $\sigma_1 = 0$ and chose several different values of $\sigma_2$.

The results based on several different values of $\sigma_2$ and 10,000 replications are reported in Table 2. For binary data, the type I error is reasonably controlled for the four statistics, even when $\sigma_2$ is large. In contrast, under the binomial distribution, random intercepts lead to discrepancies between the significance levels and the nominal level for $S_S$, $S_O$ and $LS$, while the performance of $S_P$ remains reasonable. For $S_S$ and $S_O$, the discrepancy is greater for a larger $\sigma_2$. This is because $S_O$ is suitable for testing overdispersion. It is possible that $S_S$ and $S_O$ yield large $p$-values whereas $S_P$ gives a small $p$-value. In this situation, these $p$-values suggest the presence of overdispersion instead of heterogeneity of the random effects. In addition, for the binomial distribution, $S_O$ is much more powerful than $LS$ in detecting the nuisance overdispersion. For the linear mixed model, the heterogenous variance leads to inflated type I error of $LS$, while the performance of $S_S$, $S_P$ and $S_O$ remains stable.

### 4.2. Yale family study of comorbidity of alcoholism and anxiety (YFS-CAA).

The YFSCAA was conducted to examine the patterns of familial aggregation of alcoholism in the relatives of 115 probands with alcohol dependence compared to those of 147 psychiatric (80 probands with anxiety disorders) and normal controls (67 probands with no history of psychiatric disorders). The total sample used for the familial aggregation analyses included 222 probands who had 1194 adult first-degree relatives and spouses. We refer to [23] for a detailed description of the study design and data collection. Recently, Zhang, Feng and Zhu [34] developed a latent variable model as described in Example 1 and a two-step procedure for assessing familial aggregation and heritability of disease, based on the assumption that the elements of $\mathbf{v}_i$ follow a Bernoulli distribution. The importance of our reanalysis is to demonstrate how to use our new results to remove the restrictive Bernoulli assumption on $\mathbf{v}_i$.

As in Section 2.4 of [34], for any $(i, j)$, we have

$$\mathrm{var}(b_{i,j}) = \alpha_2^2 + \gamma_1(1 - \gamma_1)[(1 - \gamma_1)\alpha_3^2 + \gamma_1(\alpha_3 + \alpha_4)^2 + (\alpha_3 + \gamma_1\alpha_4)^2].$$

Similarly, we can get $\mathrm{cov}(b_{i,j}, b_{i,k})$ for all $j, k = 1, \ldots, m_i$. Let $\sigma_T = \mathrm{var}(b_{i,j})$. Then we have $\mathrm{var}(\mathbf{b}_i) = \Sigma_{i,i} = \sigma_T W_{i,i}(\gamma)$ for the $i$th family. For example, for

TABLE 2

*Comparison of the type* I *error in the presence of nuisance overdispersion for the score test statistics under the logistic and linear mixed models at significance level* 0.05. *"Considered" and "ignored" represent including or excluding the correlation in the score test statistics*

| $\sigma_2^2$ | $LS$ | $\hat{S}_O$ | $\hat{S}_P$ | $\hat{S}_S$ | $LS$ | $\hat{S}_O$ | $\hat{S}_P$ | $\hat{S}_S$ |
|---|---|---|---|---|---|---|---|---|
| | | | | logistic mixed model (18) and Bernoulli distribution | | | | |
| | | | $n = 30$ ($m_i = 5$) | | | | $n = 50$ ($m_i = 5$) | |
| 0.3 | 0.055 | 0.054 | 0.050 | 0.054 | 0.051 | 0.056 | 0.040 | 0.040 |
| 0.6 | 0.048 | 0.054 | 0.056 | 0.066 | 0.048 | 0.054 | 0.048 | 0.046 |
| 0.9 | 0.055 | 0.056 | 0.061 | 0.074 | 0.042 | 0.046 | 0.066 | 0.048 |
| 1.2 | 0.058 | 0.058 | 0.070 | 0.068 | 0.045 | 0.060 | 0.061 | 0.060 |
| | | | | logistic mixed model (18) and binomial distribution | | | | |
| | | | $n = 30$ ($m_i = 5$) | | | | $n = 50$ ($m_i = 5$) | |
| 0.2 | 0.131 | 0.195 | 0.054 | 0.105 | 0.163 | 0.346 | 0.055 | 0.154 |
| 0.4 | 0.265 | 0.522 | 0.052 | 0.188 | 0.364 | 0.784 | 0.062 | 0.298 |
| 0.6 | 0.420 | 0.792 | 0.061 | 0.276 | 0.568 | 0.964 | 0.062 | 0.435 |
| 0.8 | 0.558 | 0.942 | 0.063 | 0.374 | 0.734 | 0.998 | 0.073 | 0.568 |
| 1.0 | 0.670 | 0.987 | 0.061 | 0.444 | 0.863 | 1.000 | 0.068 | 0.673 |
| | | | | linear mixed model (19) | | | | |
| | | | | $n = 40$ ($m_i = 4$) | | | | |
| | | | ignored | | | | considered | |
| 0.25 | 0.084 | 0.010 | 0.048 | 0.039 | 0.188 | 0.011 | 0.042 | 0.031 |
| 0.50 | 0.119 | 0.010 | 0.040 | 0.041 | 0.199 | 0.029 | 0.043 | 0.045 |
| 0.75 | 0.170 | 0.009 | 0.041 | 0.034 | 0.272 | 0.023 | 0.043 | 0.031 |
| 1.00 | 0.212 | 0.005 | 0.039 | 0.030 | 0.317 | 0.016 | 0.037 | 0.031 |

a nuclear family with two siblings, we can show that

$$\text{var}(\mathbf{b}_i) = \Sigma_{i,i} = \sigma_T \begin{pmatrix} 1 & \rho_2 & \rho_0 & \rho_0 \\ \rho_2 & 1 & \rho_0 & \rho_0 \\ \rho_0 & \rho_0 & 1 & \rho_1 \\ \rho_0 & \rho_0 & \rho_1 & 1 \end{pmatrix},$$

where $\rho_2 = \alpha_2^2/\sigma_T$, $\rho_0 = [\alpha_2^2 + \gamma_1(1-\gamma_1)(\alpha_3 + \gamma_1\alpha_4)^2]/\sigma_T$ and $\rho_1 = [\alpha_2^2 + \gamma_1(1-\gamma_1)(\alpha_3 + \gamma_1\alpha_4)^2 + 0.25\gamma_1^2(1-\gamma_1)^2\alpha_3^2]/\sigma_T$. Let $\alpha_2/S_\alpha = \cos(\gamma_2)$, $\alpha_3/S_\alpha = \sin(\gamma_2)\cos(\gamma_3)$ and $\alpha_4/S_\alpha = \sin(\gamma_2)\sin(\gamma_3)$, where $S_\alpha = \sqrt{\alpha_2^2 + \alpha_3^2 + \alpha_4^2}$. Then $\rho_0$, $\rho_1$ and $\rho_2$ can be written as functions of $\gamma_1$, $\gamma_2$ and $\gamma_3$, which are nuisance parameters here. It is noteworthy that deriving these correlation parameters is relatively straightforward for a general pedigree.

The score test statistics presented in Section 2 can be used to detect the familial correlation that includes both environmental and genetic factors through testing the hypothesis $\sigma_T = 0$. Under the null hypothesis, the maximum likelihood estimate is $(1.3341, -0.4181, 0.0178, -1.6501, -1.1522)$. To compute the maximum score test statistics, we used $r_0 = 10{,}000$ and ap-proximated the nuisance parameter domain $\Gamma = [0.01, 0.99] \times [-\pi/2, \pi/2] \times [-\pi/2, \pi/2]$ by a $15 \times 15 \times 15$ grid, $\Gamma_A = \{(i/16, j\pi/14, k\pi/14) : i = 1, \ldots, 15; j, k = -7, \ldots, 7\}$. Accordingly, $S_O$, $S_P$ and $S_S$ equal 1.23, 64.22 and 63.91, respectively. The $p$-value for $S_O$ is 0.158, revealing no evidence for the overdispersion. The $p$-values for $S_P$ and $S_S$ are less than 0.0001, providing significant support for familial aggregation and inheritability of alcoholism. To ensure that the size of the approximating grid does not affect our analysis, we considered a series of grids from smallest size $2 \times 2 \times 2$ to largest size $100 \times 100 \times 100$. The differences in the approximated values for $S_O$, $S_P$ and $S_S$ are indeed so small that they had no impact on our analysis.

**5. Discussion.** We have proposed several score statistics to test homogeneity and overdispersion in the mixed effects and latent variable models. The major advantage of these statistics is that they do not depend on the distribution of the random effects except for their mean and variance. Simulation studies demonstrate that both $S_P$ and $S_S$ have great power in detecting the heterogeneity in latent variables, but the type I error of $S_S$ is inflated in the case of nuisance overdispersion (Table 2). We have also examined a number of simulated data sets and one real application to highlight the broad spectrum of the applications for which our test procedures can be used. Another advantage of these statistics is that they automatically impose the positive semidefinite constraint on the variance components of $\Sigma_{i,i'}(\gamma)$'s. For the model in Example 2, the statistic $S_S$ reduces to the projection score test statistic [16], which asymptotically follows a mixture of $\chi^2$ distributions under a $H_0$ [26]. The simulation studies in Section 4 suggest that using a constrained score test can substantially increase the power of detecting heterogeneity. See [38] for detailed discussion.

The score statistics proposed here are to test whether all the variance components are zero. When the null hypothesis is rejected, it is also of interest to test whether some of the variance components are zero. The advantage of testing the overall hypothesis on all the variance components, followed by identifying some nonzero components, is control of the type I error. Although we do not discuss how to identify the nonzero components, our results can be useful for this purpose. For instance, $S_S$ can be directly applied to clustered designs [22] as $n \to \infty$ and when all $m_i$'s are bounded by a constant, and the asymptotic distribution of $S_S$ can also be derived under an M-dependent sequence in nested models by using the functional central limit theorem for dependent data [13]. In particular, we can follow the derivation in Section 2 to develop a score test by using the parameterization given in Section 1, and use a parametric bootstrap (or resampling procedure) to approximate the $p$-value. However, it is beyond the scope of this work to address all these related issues in detail.

Many issues still merit further research. One major issue is the empirical performance of the test statistics in finite samples under different situations, such as proportional hazard models with random effects [30] and genetic linkage test. Other possible applications include tests of spatial homogeneity for spatial processes, tests of serial correlation for state space models and tests of the Markov (or semi-Markov) hypothesis [8]. In addition, our result can be used to address important practical problems such as the selection of the random effects components in a generalized linear mixed model [7, 11]. Our results combined with those in [10] may also be useful in mixed effects models for semiparametric regression. It is noteworthy that our assumptions in this paper are not optimal. Some extensions are still possible and warrant future research. Another major issue is to further assess the impact of the grid dimension on the quality of approximation, even though some empirical evidence suggests that a rough grid works very well. Also see [1] and [36].

## APPENDIX: PROOFS AND TECHNICAL DETAILS

We introduce some notation as follows:

$$x_K = \frac{U_K}{\sqrt{EU_K^2}}, \qquad c_{K,K}(\gamma) = 0,$$

$$c_{K,K'}(\gamma) = b_{K,K'}(\gamma)\sqrt{\frac{EU_K^2 EU_{K'}^2}{I_{TP}(\gamma)}}, \qquad K \neq K',$$

$$y_K = \frac{(U_K^2 - V_K)}{\sqrt{\text{Var}(U_K^2 - V_K)}},$$

$$d_K(\gamma) = b_{K,K}(\gamma)\sqrt{\frac{\text{Var}(U_K^2 - V_K)}{I_{TO}(\gamma)}},$$

$$\lambda_{N(1)}(\gamma) = \sqrt{\frac{I_{TP}(\gamma)}{I_{TS}(\gamma)}} \quad \text{and} \quad \lambda_{N(2)}(\gamma) = \sqrt{\frac{I_{TO}(\gamma)}{I_{TS}(\gamma)}}.$$

Thus, $X_P(\gamma) = \sum_{K \neq K'} c_{K,K'}(\gamma) x_K x_{K'}$, $X_O(\gamma) = \sum_K d_K(\gamma) y_K$ and

$$X_S(\gamma) = \frac{T_S(\gamma)}{\sqrt{I_{TS}(\gamma)}} = \lambda_{N(1)}(\gamma) \sum_{K \neq K'} c_{K,K'}(\gamma) x_K x_{K'} + \lambda_{N(2)}(\gamma) \sum_K d_K(\gamma) y_K.$$

### A.1. Regularity conditions of Theorem 1.

(A1) As $N \to \infty$, $X_P(\gamma)$ converges to a Gaussian process with mean zero and covariance matrix $\rho_{(1)}(\gamma, \gamma')$ and $\lim_{N \to \infty} \sup_{\gamma \in \Gamma} |\mu_{\max}[C(\gamma)]| = 0$, where $\mu_{\max}[C(\gamma)]$ is the largest absolute eigenvalue of $C(\gamma) = (c_{K,K'}(\gamma))$.

(A2) $\lim_{N \to \infty} \sup_K \sup_{\gamma \in \Gamma} |d_K(\gamma)| = 0$. The $d_K(\gamma)$'s and $c_{K,K'}(\gamma)$'s have first-order derivatives with respect to $\gamma$ and for any $\{\gamma_t, t = 1, \ldots, q_2\}$, $\sum_K \sup_{\gamma \in \Gamma} [\partial_{\gamma_t} d_K(\gamma)]^2 < \infty$ and $\sum_{K \neq K'} \sup_{\gamma \in \Gamma} [\partial_{\gamma_t} c_{K,K'}(\gamma)]^2 < \infty$.

(A3) The sequence $\{\sup_{\gamma \in \Gamma} |d_K(\gamma)| \|y_K\| : K = 1, \ldots, N\}$ satisfies the Lindeberg condition.

(A4) For any $\gamma$ and $\gamma'$ in $\Gamma$, $\lim_{N \to \infty} \rho_{N(2)}(\gamma, \gamma') = \sum_K d_K(\gamma) d_K(\gamma') = \rho_{(2)}(\gamma, \gamma')$.

(A5) $\lim_{N \to \infty} \sup_{\gamma \in \Gamma} |\lambda_{N(1)}(\gamma) - \lambda_1(\gamma)| = 0$ and $\lambda_1(\gamma)$ is continuous in $\gamma$.

PROOF OF THEOREM 1. In terms of $\sum_K d_K(\gamma) y_K$, we can directly apply the Jain–Marcus theorem [29] by using assumptions (A2)–(A4). The finite convergence of $X_S(\gamma)$ can be observed from Theorem 5 of [15] by using assumptions (A1)–(A3). To prove the asymptotic equicontinuity of $X_S(\gamma)$, we note that for any $\gamma$ and $\gamma'$ in $\Gamma$, $|X_S(\gamma) - X_S(\gamma')|$ is bounded by

$$|\lambda_{N(1)}(\gamma) X_P(\gamma) - \lambda_{N(1)}(\gamma') X_P(\gamma')| + |\lambda_{N(2)}(\gamma) X_O(\gamma) - \lambda_{N(2)}(\gamma') X_O(\gamma')|.$$

The first term $|\lambda_{N(1)}(\gamma) X_P(\gamma) - \lambda_{N(1)}(\gamma') X_P(\gamma')|$ is bounded by

$$|\lambda_{N(1)}(\gamma) - \lambda_{N(1)}(\gamma')| |X_P(\gamma)| + |\lambda_{N(1)}(\gamma')| |X_P(\gamma) - X_P(\gamma')|.$$

From (A5), it follows that $|\lambda_{N(1)}(\gamma) - \lambda_{N(1)}(\gamma')|$ can be sufficiently small when $\gamma$ and $\gamma'$ are sufficiently close. Using the fact that $X_P(\gamma) = O_p(1)$, $\lambda_{N(1)}(\gamma') \leq 1$ and $X_P(\gamma)$ is stochastically continuous, we can prove that for any $\varepsilon, \eta > 0$, there exists a $\delta > 0$ such that

$$\lim_{N \to \infty} P\left\{ \sup_{\|\gamma - \gamma'\| \leq \delta} |\lambda_{N(1)}(\gamma) X_P(\gamma) - \lambda_{N(1)}(\gamma') X_P(\gamma')| > \varepsilon \right\} < \eta.$$

By using similar arguments, we can handle the second term. Therefore, $X_S(\gamma)$ is stochastically continuous. This completes the proof of Theorem 1. $\square$

REMARK 1. Assumption (A1) will be established in Theorem 3 where we introduce sufficient conditions (C1)–(C4). Assumptions (A2)–(A4) can be replaced by the assumptions of Theorem 10.6 of [25], but for simplicity, we prefer (A2)–(A4) because they can be easily checked for all examples considered here as well as those in [38].

Note that $U_K$ depends on $\xi_*$ implicitly. So, we denote it by $U_K(\mathbf{h} N^{-1/2})$ when we replace $\xi_*$ in $U_K$ by $\xi_* + \mathbf{h} N^{-1/2}$. We introduce similar notation for $V_K$, $\{I_{TO}(\gamma), I_{TS}(\gamma), I_{TP}(\gamma)\}$ and $\{I_{EO}(\gamma), I_{ES}(\gamma), I_{EP}(\gamma)\}$. After replacing $\xi_*$ by $\xi_* + \mathbf{h} N^{-1/2}$ in $T_P(\gamma)$, $T_S(\gamma)$ and $T_O(\gamma)$, we get $T_P(\gamma, \mathbf{h})$, $T_S(\gamma, \mathbf{h})$ and $T_O(\gamma, \mathbf{h})$. We define $\lambda_{N(3)}(\gamma) = \sqrt{I_{EP}(\gamma)/I_{ES}(\gamma)}$, $\lambda_{N(4)}(\gamma) = \sqrt{I_{EO}(\gamma)/I_{ES}(\gamma)}$, $y_K(\mathbf{h}) = [U_K(\mathbf{h}N^{-1/2})^2 - V_K(\mathbf{h}N^{-1/2})]/\sqrt{\mathrm{Var}(U_K^2 - V_K)}$,

$$x_K(\mathbf{h}) = \frac{U_K(\mathbf{h}N^{-1/2})}{\sqrt{EU_K^2}} \quad \text{and} \quad d_{O,K}(\gamma) = b_{K,K}(\gamma)\sqrt{\frac{\mathrm{Var}(U_K^2 - V_K)}{I_{EO}(\gamma)}}.$$

With these preparations, we get

$$X_P(\gamma, \mathbf{h}) = \sum_{K \neq K'} c_{K,K'}(\gamma) x_K(\mathbf{h}) x_{K'}(\mathbf{h}) \sqrt{\frac{I_{EP}(\gamma)}{I_{EP}(\gamma, \mathbf{h}N^{-1/2})}},$$

$$X_O(\gamma, \mathbf{h}) = \sum_K d_{O,K}(\gamma) y_K(\mathbf{h}) \sqrt{\frac{I_{EO}(\gamma)}{I_{EO}(\gamma, \mathbf{h}N^{-1/2})}}$$

and

$$X_S(\gamma, \mathbf{h}) = [\lambda_{N(3)}(\gamma) X_P(\gamma, \mathbf{h}) + \lambda_{N(4)}(\gamma) X_O(\gamma, \mathbf{h})] \sqrt{\frac{I_{ES}(\gamma)}{I_{ES}(\gamma, \mathbf{h}N^{-1/2})}}.$$

The following conditions are assumed for Theorem 2.

### A.2. Regularity conditions of Theorem 2.

(B1) $\lim_{N \to \infty} \sup_{\gamma \in \Gamma, \|\mathbf{h}\|_2 \leq M} \|I_{EO}(\gamma)[I_{EO}(\gamma, \mathbf{h}N^{-1/2})]^{-1} - 1\| = 0$ and

$$\lim_{N \to \infty} \sup_{\gamma \in \Gamma, \|\mathbf{h}\|_2 \leq M} \left\| \frac{I_{EP}(\gamma)}{I_{EP}(\gamma, \mathbf{h}N^{-1/2})} - 1 \right\| = 0.$$

(B2) For any $\|\mathbf{h}\|_2 \leq M$, $\sup_K \mathrm{Var}[y_K(\mathbf{h}) - y_K(\mathbf{0}) - \partial_{\mathbf{h}} y_K(0)^T \mathbf{h}] \to 0$ as $N \to \infty$ and $|y_K(\mathbf{h}) - y_K(\mathbf{h}')| \leq z_K \|\mathbf{h} - \mathbf{h}'\|_2$, where $\sup_K E(z_K^2) < \infty$. In addition, $\sup_K E[y_K(\mathbf{0})^2 + \|\partial_{\mathbf{h}} y_K(\mathbf{0})\|_2^2] < \infty$, $\sup_{\gamma \in \Gamma} \sum_{K=1}^N [d_{O,K}(\gamma)]^2 < \infty$ and $\sum_{K=1}^N \sup_{\gamma \in \Gamma} [\partial_{\gamma_t} d_{O,K}(\gamma)]^2 < \infty$.

(B3) $N^{1/2}(\hat{\xi} - \xi_*) = N^{-1/2} \sum_K F_K + o_p(1)$, $N^{-1/2} \sum_K F_K = O_p(1)$, and

$$\lim_{N \to \infty} \sup_{\gamma \in \Gamma} \left| \sum_K d_{O,K}(\gamma) \, \partial_{\mathbf{h}} y_K(\mathbf{0}) - N^{1/2} \frac{1}{\sqrt{I_{EO}(\gamma)}} J_N(\gamma) \right| = 0.$$

(B4) $f_K(\gamma) = b_{K,K}(\gamma)/\sqrt{I_{EO}(\gamma)}$ has the first-order derivative with respect to $\gamma$, and for any $\{\gamma_t, t = 1, \ldots, q_2\}$, $\sum_K \sup_{\gamma \in \Gamma} [\partial_{\gamma_t} d_K(\gamma)]^2 < \infty$. We assume similar conditions for all components of $J_{N(O)}(\gamma) = J_N(\gamma)/\sqrt{I_{EO}(\gamma)}$.

(B5) The sequence $\{\sup_{\gamma \in \Gamma} |f_K(\gamma)(U_K^2 - V_K) - J_{N(O)}(\gamma)^T F_K| : K = 1, \ldots, N\}$ satisfies the Lindeberg condition.

(B6) For any $\gamma$ and $\gamma'$ in $\Gamma$, $\lim_{N\to\infty} \rho_{N(3)}(\gamma, \gamma') = \rho_{(3)}(\gamma, \gamma')$, where $\rho_{N(3)}(\gamma, \gamma')$ is given by

$$\sum_K \operatorname{Cov}[f_K(\gamma)(U_K^2 - V_K) - J_{N(O)}(\gamma)^T F_K, f_K(\gamma')(U_K^2 - V_K) - J_{N(O)}(\gamma')^T F_K].$$

(B7) For any given $M > 0$, $\sup_{\gamma \in \Gamma, \|\mathbf{h}\|_2 \leq M} |\sum_{K \neq K'} c_{K,K'}(\gamma)[x_K(\mathbf{h})x_{K'}(\mathbf{h}) - x_K x_{K'}]| = o_p(1)$, and $X_P(\gamma, \mathbf{0})$ converges in distribution to a Gaussian process with mean zero and covariance $\rho_{(1)}(\gamma, \gamma')$.

(B8) As $N \to \infty$, $\lambda_{N(3)}(\gamma)$ uniformly converges to $\lambda_{(3)}(\gamma)$ for all $\gamma \in \Gamma$, and $\lambda_{(3)}(\gamma)$ is continuous in $\gamma$.

REMARK 2.   Some sufficient conditions for (B3) in a general mixed model have been given by Jiang [18]. Some sufficient conditions for (B7) will be given in (C1)–(C5), and other conditions will be given in Theorem 5. Also see Theorems 3 and 5 for more details.

PROOF OF THEOREM 2.   The proof of Theorem 2 consists of three steps. In the first step, we will establish that

$$(20) \qquad X_O(\gamma, \mathbf{h}) = \sum_K d_{O,K}(\gamma)[y_K(\mathbf{0}) + \partial_{\mathbf{h}} y_K(\mathbf{0})\mathbf{h}](1 + o_p(1)) + o_p(1).$$

From (B1), it follows that $X_O(\gamma, \mathbf{h}) = \sum_K d_{O,K}(\gamma)y_K(\mathbf{h})(1 + o_p(1))$. Furthermore, we consider the stochastic process $\mathrm{SP(I)} = \sum_K d_{O,K}(\gamma)[y_K(\mathbf{h}) - y_K(\mathbf{0}) - \partial_{\mathbf{h}} y_K(\mathbf{0})\mathbf{h}]$ indexed by $\{(\gamma, \mathbf{h}) : \gamma \in \Gamma, \|\mathbf{h}\|_2 \leq M\}$. For each fixed $(\gamma, \mathbf{h})$, the variance of $\mathrm{SP(I)}$ converges to zero and (B2) leads to the result that $\mathrm{SP(I)}$ is stochastically continuous. Thus, (20) is proved. We can use (B3) to deduce that $X_O(\gamma, \sqrt{N}(\hat{\xi} - \xi_*))$ can be approximated by

$$\frac{1}{\sqrt{I_{EO}(\gamma)}} \sum_K [b_{K,K}(\gamma)(U_K^2 - V_K) - J_N(\gamma)^T F_K][1 + o_p(1)] + o_p(1);$$

therefore, $\hat{X}_O(\gamma)$ converges to a Gaussian process with mean zero and covariance matrix $\rho_3(\gamma, \gamma')$ because (B4)–(B6) are sufficient conditions for this claim.

The second step is to show that

$$X_P(\gamma, \sqrt{N}(\hat{\xi} - \xi_*)) = \sum_{K \neq K'} c_{K,K'}(\gamma)x_K x_{K'}[1 + o_p(1)] + o_p(1).$$

This can be proved by using (B1) and (B7), which ends the second step. As the last step, we combine the results on $X_P(\gamma, \mathbf{h})$ and $X_O(\gamma, \mathbf{h})$ and then follow the proof of Theorem 1 to complete the proof for Theorem 2.   □

**A.3. Regularity conditions of Theorem 3.**

(C1) Let $x_1, \ldots, x_N$ be a sequence of independent random variables such that $Ex_K = 0$ and $Ex_K^2 = 1$ for all $K = 1, \ldots, N$. We assume that $\sup_K \{E|x_K|^p\} < \infty$ for some integer $p$ greater than $\max(q_2, 2)$.

(C2) $\mathrm{Var}[X_P(\gamma)] = 2\sum_{K \neq K'} c_{K,K'}^2(\gamma) = 2$ and $\lim_{N \to \infty} \sup_{\gamma \in \Gamma} |\mu_{\max}[C(\gamma)]| = 0$.

(C3) $c_{K,K'}(\gamma)$ has continuous first-order and second-order derivatives with respect to $\gamma$. Let $\gamma_t$ be any component of $\gamma$. We assume that $\sum_{K \neq K'} \sup_{\gamma \in \Gamma} [\partial_{\gamma_t} c_{K,K'}(\gamma)]^2 < \infty$ for $t = 1, \ldots, q_2$.

(C4) For any $\gamma$ and $\gamma'$ in $\Gamma$,

$$\lim_{N \to \infty} \rho_{N(1)}(\gamma, \gamma') = \sum_{K \neq K'} c_{K,K'}(\gamma) c_{K,K'}(\gamma') = \rho_{(1)}(\gamma, \gamma').$$

PROOF OF THEOREM 3. First, we need to show that any finite-dimensional distributions of $\{X_P(\gamma) : \gamma \in \Gamma\}$ converge weakly to the corresponding finite-dimensional distributions of $\{G_P(\gamma) : \gamma \in \Gamma\}$. From (C1)–(C3) and the martingale convergence theorem, it follows that $\sum_{K \neq K'} c_{K,K'}(\gamma) x_K x_{K'}$ converges to the standard normal in distribution for any $\gamma \in \Gamma$; see [24] and [15].

Let us consider two points $\gamma_1$ and $\gamma_2$ in $\Gamma$. By using the Cramér–Wald device, we need to show that for any $a_1$ and $a_2$ in $R$,

$$X_P(\gamma_1, \gamma_2) = a_1 X_P(\gamma_1) + a_2 X_P(\gamma_2) \xrightarrow{L} N[0, 2a_1^2 + 2a_2^2 + 4a_1 a_2 \rho_{(1)}(\gamma_1, \gamma_2)].$$

From (C2) and (C4), we know that $\mathrm{Var}[X_P(\gamma_1, \gamma_2)]$ converges to $2a_1^2 + 2a_2^2 + 4a_1 a_2 \rho_{(1)}(\gamma_1, \gamma_2)$. If $a_1^2 + a_2^2 + 2a_1 a_2 \rho_{(1)}(\gamma_1, \gamma_2) = 0$, then $X_P(\gamma_1, \gamma_2)$ converges to zero in probability and $a_1 G_P(\gamma_1) + a_2 G_P(\gamma_2) = 0$. In other cases, we have $a_1^2 + a_2^2 + 2a_1 a_2 \rho_{(1)}(\gamma_1, \gamma_2) > 0$. From (C2), it follows that

$$|\mu_{\max}[a_1 C(\gamma_1) + a_2 C(\gamma_2)]| \leq |a_1| |\mu_{\max}[C(\gamma_1)]| + |a_2| |\mu_{\max}[C(\gamma_2)]| \to 0.$$

Thus, $X_P(\gamma_1, \gamma_2)$ converges to the desired normal random variable in distribution. Similarly, we can generalize this result to any finite cases.

From Lemma 1.3 of [24], it follows that $\{E|\sum_{K \neq K'} [c_{K,K'}(\gamma_1) - c_{K,K'}(\gamma_2)] \times x_K x_{K'}|^p\}^{2/p} \leq C \sum_{K \neq K'} [c_{K,K'}(\gamma_1) - c_{K,K'}(\gamma_2)]^2$, where $C$ is a scalar independent of $N$. By using (C3), we have $\{E|\sum_{K \neq K'} [c_{K,K'}(\gamma_1) - c_{K,K'}(\gamma_2)] \times x_K x_{K'}|^p\}^{1/p} \leq C \|\gamma_1 - \gamma_2\|_2$. To prove the stochastic equicontinuity of $X_P(\gamma)$, we just need to show that $E \sup_{\|\gamma_1 - \gamma_2\|_2 \leq \delta} |\sum_{K \neq K'} [c_{K,K'}(\gamma_1) - c_{K,K'}(\gamma_2)] \times x_K x_{K'}|^p \to 0$ as $\delta \to 0$ and $N \to \infty$. We can finish our proof by noting that $\Gamma$ is a bounded compact set of $R^{q_2}$, whose packing number $D(t, \Gamma, \|\cdot\|_2)$ is of the order of $t^{-q_2}$. Theorem 2.2.4 of [29] concludes the proof. $\square$

PROOF OF THEOREM 4. First, we will prove the unconditional weak convergence of $X_P^{(r)}(\gamma)$. After some calculations, we can show that $E[X_P^{(r)}(\gamma)] =$ 0 and $\mathrm{Cov}[X_P^{(r)}(\gamma_1), X_P^{(r)}(\gamma_2)] = 2\sum_{K\neq K'} c_{K,K'}(\gamma_1)c_{K,K'}(\gamma_2) = 2\rho_{N(1)}(\gamma_1,\gamma_2)$. Therefore, $X_P^{(r)}(\gamma)$ and $X_P(\gamma)$ have the same mean and covariance structures. Following the proof of Theorem 4 in [15], we can show that $X_P^{(r)}(\gamma)$ as the sum of martingale differences converges to a normal random variable for any $\gamma \in \Gamma$, as $N \to \infty$. The Cramér–Wald device is applicable in any finite case. Following the similar argument in Theorem 3, we can show the stochastic continuity of $X_P^{(r)}(\gamma)$. Therefore, $X_P^{(r)}(\gamma)$ converges to $G_P(\gamma)$ in distribution; that is, $X_P^{(r)}(\gamma)$ is asymptotically measurable.

Second, given $x_1,\ldots,x_N$, we have $X_P^{(r)}(\gamma) \sim N[0, 2\sum c_{K,K'}(\gamma)^2 x_K^2 x_{K'}^2]$ and $E_V[X_P^{(r)}(\gamma)X_P^{(r)}(\gamma')]$ is equal to $2\sum_{K\neq K'} c_{K,K'}(\gamma)c_{K,K'}(\gamma')x_K^2 x_{K'}^2$. We write $E_V[X_P^{(r)}(\gamma)X_P^{(r)}(\gamma')]/2$ as the sum of $\sum_{K\neq K'} c_{K,K'}(\gamma)c_{K,K'}(\gamma')(x_K^2-1)(x_{K'}^2-1)$, $2\sum_{K\neq K'} c_{K,K'}(\gamma)c_{K,K'}(\gamma')(x_{K'}^2-1)$ and $\rho_{N(1)}(\gamma,\gamma')$. The first term is also a random quadratic form. Its mean is zero and its variance is bounded by $C\max_K\{\sum_{K'=1}^N c_{K,K'}(\gamma')^2\}$, which converges to zero; see Lemma 1.2 of [24]. By using Theorem 1 of [31], we can show that $\sum_{K\neq K'} c_{K,K'}(\gamma)c_{K,K'}(\gamma')(x_K^2-1)(x_{K'}^2-1)$ converges to zero in probability. The same technique can be used to show that $\sum_{K\neq K'} c_{K,K'}(\gamma)c_{K,K'}(\gamma')(x_K^2-1)$ converges to zero in probability. Thus, $E_V[X_P^{(r)}(\gamma)X_P^{(r)}(\gamma')] \to \rho_{(1)}(\gamma,\gamma')$ in probability. We can obtain the marginal convergence in the conditional central limit theorem by using the Cramér–Wald method.

For each $\delta > 0$, let $\Gamma_\delta$ assign to each $\gamma \in \Gamma$ a closest element of a given finite $\delta$-net of $\Gamma$ with respect to $\|\cdot\|_2$. The above finite convergence results lead to $\sup_{h\in BL_1(\ell^\infty(\Gamma))} |E_V h(X_P^{(r)}(\Gamma_\delta(\cdot))) - Eh(G_P(\Gamma_\delta(\cdot)))| \to 0$ in probability, as $N \to 0$. By continuity of $G_P(\gamma)$, we have $G_P(\Gamma_\delta(\gamma)) \to G_P(\gamma)$ almost surely, as $\delta \to 0$; that is, $\lim_{\delta\to 0}\sup_{h\in BL_1(\ell^\infty(\Gamma))} |Eh(G_P(\Gamma_\delta(\cdot))) - Eh(G_P(\cdot))| = 0$. Finally, $\sup_{h\in BL_1(\ell^\infty(\Gamma))} |E_V h(X_P^{(r)}(\Gamma_\delta(\cdot))) - E_V h(X_P^{(r)}(\cdot))|$ is bounded by $E_V(\sup_{\|\gamma-\gamma'\|_2\leq\delta} |X_P^{(r)}(\gamma') - X_P^{(r)}(\gamma)|)$. Because the expectation on the left-hand side is smaller than $E(\sup_{\gamma,\gamma'\in\Gamma; \|\gamma-\gamma'\|_2\leq\delta} |X_P^{(r)}(\gamma') - X_P^{(r)}(\gamma)|)$, which was established by the unconditional weak convergence of $X_P^{(r)}(\gamma)$, the desired results follow. □

Next, we state a few more assumptions. Let $\tilde{U}_K(\mathbf{h}N^{-1/2}) = U_K(s_K, \mathbf{h}N^{-1/2}) - \mu_K(\mathbf{h}N^{-1/2})$, where $\mu_K(\mathbf{h}N^{-1/2}) = EU_K(s_K, \mathbf{h}N^{-1/2})$.

### A.4. Regularity conditions of Theorem 5.

(C5) $U_K(s_K, \xi)$ has continuous first-order and second-order derivatives with respect to $\xi$ in an open neighborhood of $\xi_*$, denoted by $\partial_\xi U_K(s_K, \xi)$ and $\partial_\xi^2 U_K(s_K, \xi)$, respectively.

(C6) $\sup_K \mathrm{Var}\{\tilde{U}_K(\mathbf{h}N^{-1/2})\tilde{U}_{K'}(\mathbf{h}N^{-1/2}) - x_K x_{K'}\} \to 0$ as $N \to \infty$.

(C7) $\sup_{\|\mathbf{h}\|_2 \le M} \sum_K \mu_K(\mathbf{h}N^{-1/2})^2 < \infty$, $\sup_K \sup_{\|\mathbf{h}\|_2 \le M} E^{1/p}|\tilde{U}_K(\mathbf{h} \times N^{-1/2})|^p < \infty$ and $\sup_K E^{1/p}|\tilde{U}_K(\mathbf{h}N^{-1/2}) - \tilde{U}_K(\mathbf{h}'N^{-1/2})|^p < c\|\mathbf{h}' - \mathbf{h}\|_2$ for some integer $p > q_2 + q_3$.

(C8) $\sup_{\gamma \in \Gamma, \|\mathbf{h}\|_2 \le M} |\sum_{K \ne K'} c_{K,K'}(\gamma)\mu_K(\mathbf{h}N^{-1/2})\tilde{U}_K(\mathbf{h}N^{-1/2})| = o_p(1)$.

PROOF OF THEOREM 5. We see that $X_P(\gamma, \mathbf{h}N^{-1/2})$ can be written as the sum of $\sum_{K \ne K'} c_{K,K'}(\gamma)\tilde{U}_K(\mathbf{h}N^{-1/2})\tilde{U}_{K'}(\mathbf{h}N^{-1/2})$,

$$2 \sum_{K \ne K'} c_{K,K'}(\gamma)\mu_K(\mathbf{h}N^{-1/2})\tilde{U}_{K'}(\mathbf{h}N^{-1/2})$$

and $\sum_{K \ne K'} c_{K,K'}(\gamma)\mu_K(\mathbf{h}N^{-1/2})\mu_{K'}(\mathbf{h}N^{-1/2})$. In the following, we will prove that every term in the foregoing equation converges to zero in probability.

For the third term, we have

$$\text{term (III)} \le \sup_{\gamma \in \Gamma} \mu_{\max}[C(\gamma)] \sup_{\|\mathbf{h}\|_2 \le M} \sum_K \mu_K(\mathbf{h}N^{-1/2})^2,$$

which converges to zero as $N$ is sufficiently large. The second term (II) is just assumption (C8).

For the first term, we need to consider the process $T_N(\gamma, \mathbf{h}) = \text{term (I)} - X_P(\gamma, \mathbf{0})$. For each $\gamma$ and $\mathbf{h}$, $T_N(\gamma, \mathbf{h})$ has mean zero and variance given by

$$2 \sum_{K \ne K'} c_{K,K'}(\gamma)^2 \mathrm{Var}\{\tilde{U}_K(\mathbf{h}N^{-1/2})\tilde{U}_{K'}(\mathbf{h}N^{-1/2}) - x_K x_{K'}\},$$

which converges to zero by assumption (C6). To establish stochastic continuity of $T_N(\gamma, \mathbf{h})$, we find that $T_N(\gamma, \mathbf{h}) - T_N(\gamma', \mathbf{h}') = (a) + (b) + (c) + (d)$, where each term on the right-hand side is given by

$$(a) = \sum_{K \ne K'} [c_{K,K'}(\gamma) - c_{K,K'}(\gamma')]\tilde{U}_K(\mathbf{h}N^{-1/2})[\tilde{U}_{K'}(\mathbf{h}N^{-1/2}) - x_{K'}],$$

$$(b) = \sum_{K \ne K'} [c_{K,K'}(\gamma) - c_{K,K'}(\gamma')]x_{K'}[\tilde{U}_K(\mathbf{h}N^{-1/2}) - x_K],$$

$$(c) = \sum_{K \ne K'} c_{K,K'}(\gamma')\tilde{U}_K(\mathbf{h}'N^{-1/2})[\tilde{U}_{K'}(\mathbf{h}N^{-1/2}) - \tilde{U}_{K'}(\mathbf{h}'N^{-1/2})],$$

$$(d) = \sum_{K \ne K'} c_{K,K'}(\gamma')\tilde{U}_{K'}(\mathbf{h}N^{-1/2})[\tilde{U}_K(\mathbf{h}'N^{-1/2}) - \tilde{U}_K(\mathbf{h}N^{-1/2})].$$

Using the same technique as in Lemma 1.3 of [24], we can finishes the proof by using Theorem 2.2.4 of [29] and assumption (C7). □

**Acknowledgments.** We thank the Editor, an Associate Editor and three anonymous referees for valuable suggestions which helped improve our presentation greatly. We thank Dr. Kathleen Merikangas for making her alcoholism data set available to us. Reprints can be requested via e-mail: heping.zhang@yale.edu.

MRI Unit
Department of Psychiatry
Columbia University Medical Center
and
New York State Psychiatric Institute
1051 Riverside Drive
New York, New York 10032
USA
E-mail: hz2114@columbia.edu

Department of Epidemiology
    and Public Health
Yale University School of Medicine
60 College Street
New Haven, Connecticut 06520-8034
USA
and
Jiangxi Normal University
Nanchang
China
E-mail: heping.zhang@yale.edu